\documentclass[12pt]{article}
\usepackage{amsmath}
\usepackage{amssymb}
\usepackage{array}
\usepackage{geometry}
\usepackage{graphicx}
\usepackage{enumerate}
\usepackage[small, compact] {titlesec}
\numberwithin{equation}{section}
\begin{document}
\vskip7cm\noindent
\begin{center}{\bf Explicit Evaluations of Matrix-variate Gamma}\\
\vskip.2cm{\bf and Beta Integrals in the Real and Complex Cases}\\
\vskip.3cm{A.M. Mathai}\\
 \vskip.2cm{Director, Centre for
Mathematical Sciences India}\\
 \vskip.1cm{[Arunapuram
P.O., Palai, Kerala-686574, Kerala, India]}\\
directorcms458@gmail.com , Phone:91+9495427558\\
\vskip.2cm{and}\\
 \vskip.2cm{Emeritus Professor of
Mathematics and Statistics, McGill University Canada;\\
mathai@math.mcgill.ca}\\
 \vskip.1cm{[805 Sherbrooke
Street West, Montreal, Quebec, Canada, H3A2K6]}
\end{center}
\newpage
\vskip.5cm\noindent{\bf Abstract} \vskip.3cm Matrix transformations in terms of triangular matrices is the easiest method of evaluating
matrix-variate gamma and beta integrals in the real and complex cases. Here we give several procedures of explicit evaluation of gamma and beta integrals in the general real and complex situations. The procedure also
reveals the structure of these matrix-variate integrals. Apart from the evaluation of matrix-variate gamma and beta integrals, the procedure can also be applied to evaluate such integrals explicitly in similar situations. Various methods described here will be useful to those who are working on integrals involving real-valued scalar functions of matrix argument in general and gamma and beta integrals in particular.


\newpage
\vskip.5cm\noindent{\bf 1.\hskip.3cm Introduction}
\vskip.3cm First we consider matrix-variate gamma integrals in the real case, then we look at matrix-variate type-1 beta integrals in the real case. The procedure is parallel in the matrix-variate type-2 beta integrals. Then we look at all these in the complex domain.

\vskip.3cm\noindent{\bf 1.1.\hskip.3cm Real matrix variate gamma integral}
\vskip.3cm
Matrix-variate gamma integral is a very popular integral in many areas. A particular case is the most popular Wishart density in multivariate statistical analysis. Let $X$ be a $p\times p$ real symmetric and positive definite matrix of mathematical or random variables. Consider the real-valued function of matrix argument
$$f(X)=C~|X|^{\alpha-\frac{p+1}{2}}{\rm e}^{-{\rm tr}(BX)}\eqno(1.1)
$$where $C$ is a constant, $|(\cdot)|$ denotes the determinant of $(\cdot)$ and ${\rm tr}(\cdot)$ denotes the trace of the matrix $(\cdot)$. All matrices appearing in this article are $p\times p$ unless stated otherwise. When $X$ is real and positive definite, $X>O$, then $f(X)$ in (1.1) represents a real matrix-variate gamma density when $C=\frac{|B|^{\alpha}}{\Gamma_p(\alpha)}$ where $B>O$ is a constant matrix and
$$\Gamma_p(\alpha)=\pi^{\frac{p(p-1)}{4}}\Gamma(\alpha)\Gamma(\alpha-\frac{1}{2})...\Gamma(\alpha-\frac{p-1}{2}),
~\Re(\alpha)>\frac{p-1}{2}.\eqno(1.2)
$$When $B$ is of the form $B={\frac{1}{2}}V^{-1}, ~V=V'>O$, where a prime denotes the transpose, then $f(X)$ in (1.1) is the Wishart density in multivariate statistical analysis, which is the central density in the area, see for example, Anderson [1], Kshirsagar [2], Srivastava and Khatri [3]. The real  matrix-variate gamma integral is
$$\Gamma_p(\alpha)=\int_{X>O}|X|^{\alpha-\frac{p+1}{2}}{\rm e}^{-{\rm tr}(X)}{\rm d}X\eqno(1.3)
$$where ${\rm d}X$ is the wedge product of the $p(p+1)/2$ differentials ${\rm d}X=\prod_{i\ge j}\wedge{\rm d}x_{ij}$. For evaluating the integral in (1.3), the standard technique used is to write $X=TT'$ where $T$ is a lower or upper triangular matrix with  positive diagonal elements. Then the integral on the right of (1.3) will split into conventional integrals on individual scalar variables. When $T$ is lower triangular then the integral over $t_{ij},i>j$ gives $\sqrt{\pi}$ and there are $p(p-1)/2$ such factors giving $\pi^{\frac{p(p-1)}{4}}$. The integral over $t_{jj}>0$ gives $\Gamma(\alpha-\frac{j-1}{2}),j=1,...,p$ and the product thus gives $\Gamma_p(\alpha)$ on the left of (1.3).
 \vskip.2cm When Wishart density is derived, starting from samples from a Gaussian population, the basic technique is the triangularization process. Can we evaluate the integral on the right of (1.3) by using conventional methods, or by direct evaluation? We will look into this problem by using the technique of partitioned matrices. Let us partition
 $$X=\left[ \begin{matrix}X_{11}&X_{12}\\
 X_{21}&X_{22}
 \end{matrix}\right]
 $$where let $X_{22}=x_{pp}$ so that $X_{21}=(x_{p1},...,x_{pp-1}), X_{12}=X_{21}'$. Then
 $$|X|^{\alpha-\frac{p+1}{2}}=|X_{11}|^{\alpha-\frac{p+1}{2}}[x_{pp}-X_{21}X_{11}^{-1}X_{12}]^{\alpha-\frac{p+1}{2}}
 $$by using partitioned matrix and determinant. Note that when $X$ is positive definite, that is, $X>O$, then $X_{11}>O,x_{pp}>0$ and the quadratic form $X_{21}X_{11}^{-1}X_{12}>0$. Note that
 $$[x_{pp}-X_{21}X_{11}^{-1}X_{12}]^{\alpha-\frac{p+1}{2}}=x_{pp}^{\alpha-\frac{p+1}{2}}[1-x_{pp}^{-\frac{1}{2}}
 X_{21}X_{11}^{-\frac{1}{2}}X_{11}^{-\frac{1}{2}}X_{12}x_{pp}^{-\frac{1}{2}}]^{\alpha-\frac{p+1}{2}}.
$$Let $Y=x_{pp}^{-\frac{1}{2}}X_{21}X_{11}^{-\frac{1}{2}}$ then ${\rm d}Y=x_{pp}^{-\frac{p-1}{2}}|X_{11}|^{-\frac{1}{2}}{\rm d}X_{21}$ for fixed $X_{11},x_{pp}$, see Mathai ([4], Theorem 1.18.) The integral over $x_{pp}$ gives
$$\int_0^{\infty}x_{pp}^{\alpha+\frac{p-1}{2}-\frac{p+1}{2}}{\rm e}^{-x_{pp}}{\rm d}x_{pp}=\Gamma(\alpha),~\Re(\alpha)>0.
$$Let $u=YY'$. Then from Theorem 2.16 and Remark 2.13 of [4] and after integrating out over the Stiefel manifold we have
$${\rm d}Y=\frac{\pi^{\frac{p-1}{2}}}{\Gamma(\frac{p-1}{2})}u^{\frac{p-1}{2}-1}{\rm d}u.
$$(Note that $n$ in Theorem 2.16 corresponds to $p-1$ and $p$ corresponds to $1$). Then the integral over $u$ gives
$$\int_0^1u^{\frac{p-1}{2}-1}(1-u)^{\alpha-\frac{p+1}{2}}{\rm d}u=\frac{\Gamma(\frac{p-1}{2})\Gamma(\alpha-\frac{p-1}{2})}{\Gamma(\alpha)},~\Re(\alpha)>\frac{p-1}{2}.
$$Now, collecting all the factors, we have
\begin{align}
|X_{11}|^{\alpha+\frac{1}{2}-\frac{p+1}{2}}\Gamma(\alpha)&\frac{\pi^{\frac{p-1}{2}}}{\Gamma(\frac{p-1}{2})}
\frac{\Gamma(\frac{p-1}{2})\Gamma(\alpha-\frac{p-1}{2})}{\Gamma(\alpha)}\nonumber\\
&=|X_{11}^{(1)}|^{\alpha+\frac{1}{2}-\frac{p+1}{2}}
\pi^{\frac{p-1}{2}}\Gamma(\alpha-\frac{p-1}{2})\nonumber
\end{align}for $\Re(\alpha)>\frac{p-1}{2}$. Note that $|X_{11}^{(1)}|$ is $(p-1)\times (p-1)$ and $|X_{11}|$ after the completion of the first part of the operations is denoted by $|X_{11}^{(1)}|$, and the exponent is changed to $\alpha+\frac{1}{2}-\frac{p+1}{2}$. Now repeat the process by separating $x_{p-1,p-1}$, that is by writing
$$X_{11}^{(1)}=\left[\begin{matrix}X_{11}^{(2)}&X_{12}^{(2)}\\
X_{21}^{(2)}&x_{p-1,p-1}
\end{matrix}\right].
$$Here $X_{11}^{(2)}$ is of order $(p-2)\times (p-2)$ and $X_{21}^{(2)}$ is of order $1\times (p-2)$. As before, let $u=YY', Y=x_{p-1,p-1}^{-\frac{1}{2}}X_{21}^{(2)}[X_{11}^{(2)}]^{-\frac{1}{2}}.$ Then ${\rm d}Y=x_{p-1,p-1}^{-\frac{p-2}{2}}|X_{11}^{(2)}|^{-\frac{1}{2}}{\rm d}X_{21}^{(2)}.$ Integral over the Stiefel manifold gives $\frac{\pi^{\frac{p-2}{2}}}{\Gamma(\frac{p-2}{2})}u^{\frac{p-2}{2}-1}{\rm d}u$ and the factor containing $(1-u)$ is $(1-u)^{\alpha+\frac{1}{2}-\frac{p+1}{2}}$ and the integral over $u$ gives
$$\int_0^1u^{\frac{p-2}{2}-1}(1-u)^{\alpha+\frac{1}{2}-\frac{p+1}{2}}{\rm d}u=\frac{\Gamma(\frac{p-2}{2})\Gamma(\alpha-\frac{p-2}{2})}{\Gamma(\alpha)}.
$$Intgral over $v=x_{p-1,p-1}$ gives
$$\int_0^1v^{\alpha+\frac{1}{2}+\frac{p-2}{2}-\frac{p+1}{2}}{\rm e}^{-v}{\rm d}v=\Gamma(\alpha),~\Re(\alpha)>0.
$$Taking all product we have
$$|X_{11}^{(2)}|^{\alpha+1-\frac{p+1}{2}}\pi^{\frac{p-2}{2}}\Gamma(\alpha-\frac{p-2}{2}),~\Re(\alpha)>\frac{p-2}{2}.
$$Successive evaluations by using the same procedure gives the exponent of $\pi$ as $\frac{p-1}{2}+\frac{p-2}{2}+...+\frac{1}{2}=\frac{p(p-1)}{4}$ and the gamma product is $\Gamma(\alpha-\frac{p-1}{2})\Gamma(\alpha-\frac{p-2}{2})...\Gamma(\alpha)$ and the final result is $\Gamma_p(\alpha)$. Hence the result is verified.

\vskip.3cm\noindent{\bf 1.2.\hskip.3cm Evaluation of matrix-variate gamma in the complex case}

\vskip.3cm In the complex case, the matrices and gamma will be denoted with a tilde. In the complex case, all matrices appearing in the integrals will be $p\times p$ hermitian positive definite unless stated otherwise, denoted by $\tilde{X}>O$. Our integral of interest is
$$\tilde{\Gamma_p}(\alpha)=\int_{\tilde{X}>O}|{\rm det}(\tilde{X})|^{\alpha-p}{\rm e}^{-{\rm tr}(\tilde{X})}{\rm d}\tilde{X}.\eqno(1.4)
$$One standard procedure to evaluate the integral in (1.4) is to write the hermitian positive definite matrix as
$\tilde{X}=\tilde{T}\tilde{T}^{*}$ where $\tilde{T}$ is a lower triangular matrix with real and positive diagonal elements $t_{jj}>0,j=1,...,p$, where * indicates the conjugate transpose. Then the Jacobian can be seen to be the following, see also ([4], Theorem 3.7):
$${\rm d}\tilde{X}=2^p\{\prod_{j=1}^pt_{jj}^{2(p-j)+1}\}{\rm d}\tilde{T}\eqno(1.5)
$$and then
\begin{align}
{\rm tr}(\tilde{X})&={\rm tr}(\tilde{T}\tilde{T}^{*})\nonumber\\
&=t_{11}^2+...+t_{pp}^2+|\tilde{t_{21}}|^2+...+|\tilde{t_{p1}}|^2+...+|\tilde{t_{pp-1}}|^2\nonumber
\end{align}and
$${\rm tr}(\tilde{X}){\rm d}\tilde{X}=2^p\{\prod_{j=1}^pt_{jj}^{2\alpha-2j+1}\}{\rm d}\tilde{T}.
$$Now, integrating out over $\tilde{t_{jk}}$ for $j>k$
$$\int_{\tilde{t_{jk}}}{\rm e}^{-|\tilde{t_{jk}}|^2}{\rm d}\tilde{t_{jk}}=\int_{-\infty}^{\infty}\int_{-\infty}^{\infty}{\rm e}^{-(t_{jk1}^2+t_{jk2}^2)}{\rm d}t_{jk1}\wedge{\rm d}t_{jk2}=\pi
$$and

$$\prod_{j>k}\pi=\pi^{\frac{p(p-1)}{2}}.
$$Now,
$$2\int_0^{\infty}t_{jj}^{2\alpha-2j+1}{\rm e}^{-t_{jj}^2}{\rm d}t_{jj}=\Gamma(\alpha-j+1),~\Re(\alpha)>j-1,
$$for $j=1,...,p.$
 Now the product of all these gives
$$\pi^{\frac{p(p-1)}{2}}\Gamma(\alpha)\Gamma(\alpha-1)...\Gamma(\alpha-p+1)=\tilde{\Gamma_p}(\alpha),~\Re(\alpha)>p-1
$$and hence the result is verified.

\vskip.3cm\noindent{\bf 1.3.\hskip.3cm An alternate method based on partitioned matrix}
\vskip.3cm Let us separate $x_{pp}$. When $\tilde{X}$ is $p\times p$ hermitian positive definite then all its diagonal elements are real and positive. That is, $x_{jj}>0,j=1,...,p$. Let
$$\tilde{X}=\left[\begin{matrix}\tilde{X_{11}}&\tilde{X_{12}}\\
\tilde{X_{21}}&x_{pp}
\end{matrix}\right]
$$where $\tilde{X_{11}}$ is $(p-1)\times (p-1)$ and
$$|{\rm det}(\tilde{X})|^{\alpha-p}=|{\rm det}(\tilde{X_{11}})|^{\alpha-p}|x_{pp}-\tilde{X_{21}}\tilde{X_{11}}^{-1}\tilde{X_{12}}|^{\alpha-p}
$$and
$${\rm tr}(\tilde{X})={\rm tr}(\tilde{X_{11}})+x_{pp}.
$$Then

$$|x_{pp}-\tilde{X_{21}}\tilde{X_{11}}^{-1}\tilde{X_{12}}|^{\alpha-p}=x_{pp}^{\alpha-p}|1-x_{pp}^{-\frac{1}{2}}
\tilde{X_{21}}\tilde{X_{11}}^{-\frac{1}{2}}\tilde{X_{11}}^{-\frac{1}{2}}\tilde{X_{12}}x_{pp}^{-\frac{1}{2}}|^{\alpha-p}.
$$Put

$$\tilde{Y}=x_{pp}^{-\frac{1}{2}}\tilde{X_{21}}\tilde{X_{11}}^{-\frac{1}{2}}\Rightarrow {\rm d}\tilde{Y}=x_{pp}^{-(p-1)}|{\rm det}(\tilde{X_{11}})|^{-1}{\rm d}\tilde{X_{21}}
$$from ([4], Theorem 3.2(c)).  Now, the integral over $x_{pp}$ gives

$$\int_0^{\infty}x_{pp}^{\alpha-p+(p-1)}{\rm e}^{-x_{pp}}{\rm d}x_{pp}=\Gamma(\alpha),~\Re(\alpha)>0.
$$Let $u=\tilde{Y}\tilde{Y}^{*}$. Then ${\rm d}\tilde{Y}=u^{p-2}\frac{\pi^{p-1}}{\Gamma(p-1)}{\rm d}u$ by using Corollaries 4.5.2 and 4.5.3 of [4]. Note that $u$ is real and positive. Integral over $u$ gives

$$\int_0^{\infty}u^{(p-1)-1}(1-u)^{\alpha-(p-1)-1}{\rm d}u=\frac{\Gamma(p-1)\Gamma(\alpha-(p-1))}{\Gamma(\alpha)},\Re(\alpha)>p-1.
$$Taking the product we have
\begin{align}
|{\rm det}(\tilde{X}_{11}^{(1)})|^{\alpha+\frac{1}{2}-p}\Gamma(\alpha)&\frac{\pi^{p-1}}{\Gamma(p-1)}\frac{\Gamma(p-1)
\Gamma(\alpha-(p-1))}{\Gamma(\alpha)}\nonumber\\
&=\pi^{p-1}\Gamma(\alpha-(p-1))|{\rm det}(\tilde{X}_{11}^{(1)})|^{\alpha+\frac{1}{2}-p}\nonumber
\end{align}
where $\tilde{X}_{11}^{(1)}$ indicates $\tilde{X}_{11}$ after the first set of integrations. Now for the second stage, separate $x_{p-1,p-1}$ and the first $(p-2)\times (p-2)$ block may be denoted by $\tilde{X}_{11}^{(2)}$. Now proceed as before to get $|{\rm det}(\tilde{X}_{11}^{(2)})|^{\alpha+1-p}\pi^{p-2}\Gamma(\alpha-(p-2))$. Proceeding like this we have the exponent of $\pi$ as $(p-1)+(p-2)+...+1=p(p-1)/2$ and the gamma product will be $\Gamma(\alpha-(p-1))\Gamma(\alpha-(p-2))...\Gamma(\alpha)$ for $\Re(\alpha)>p-1$. That is,
$$\pi^{\frac{p(p-1)}{2}}\Gamma(\alpha)\Gamma(\alpha-1)...\Gamma(\alpha-(p-1))=\tilde{\Gamma_p}(\alpha).$$
\vskip.3cm\noindent{\bf 2.\hskip.3cm Evaluation of Matrix-variate Beta Integrals}
\vskip.3cm Here we will consider a direct way of evaluating matrix-variate type-1 and type-2 beta integrals in the real and complex cases.

\vskip.3cm\noindent{\bf 2.1.\hskip.3cm Evaluation of matrix-variate type-1 beta integral in the real case}

\vskip.3cm The real matrix-variate type-1 beta function is denoted by
$$B_p(\alpha,\beta)=\frac{\Gamma_p(\alpha)\Gamma_p(\beta)}{\Gamma_p(\alpha+\beta)},~\Re(\alpha)>p-1,~\Re(\beta)>p-1
$$and it has the following type-1 beta integral representation:
$$B_p(\alpha,\beta)=\int_{O<X<I}|X|^{\alpha-\frac{p+1}{2}}|I-X|^{\beta-\frac{p+1}{2}}{\rm d}X,
$$for $\Re(\alpha)>\frac{p-1}{2},~\Re(\beta)>\frac{p-1}{2}$
 where $X$ is real symmetric and positive definite $p\times p$ matrix. The standard derivation of this integral is from the properties of real matrix-variate gamma integrals by making suitable transformations, see for example, [4]. Is it possible to evaluate the integral directly and show that it is equal to $\frac{\Gamma_p(\alpha)\Gamma_p(\beta)}{\Gamma_p(\alpha+\beta)}$, where, for example,
$$\Gamma_p(\alpha)=\pi^{\frac{p(p-1)}{4}}\Gamma(\alpha)\Gamma(\alpha-\frac{1}{2})...\Gamma(\alpha-\frac{p-1}{2}),
~\Re(\alpha)>\frac{p-1}{2}?
$$For evaluating real matrix-variate gamma integral an easy method is to make  the transformation $X=TT'$ where $T$ is a lower triangular matrix with positive diagonal  elements. Even if this transformation is applied here, the integral does not simplify due to the presence of the factor $|I-X|^{\beta-\frac{p+1}{2}}$. Hence we will try to evaluate the integral by using a partitioning of the matrices and then integrating step by step. Let $X=(x_{ij})$ be a $p\times p$ matrix. Let us separate $x_{pp}$. This can be done by partitioning $|X|$ and $|I-X|$. That is, let
$$X=\left[\begin{matrix}X_{11}&X_{12}\\
X_{21}&X_{22}\end{matrix}\right]
$$where $X_{11}$ is the $(p-1)\times (p-1)$ leading submatrix, $X_{21}$ is $1\times (p-1)$, $X_{22}=x_{pp}$ and $X_{12}=X_{21}'$. Then $|X|=|X_{11}|[x_{pp}-X_{21}X_{11}^{-1}X_{12}]$ and
$$
|X|^{\alpha-\frac{p+1}{2}}=|X_{11}|^{\alpha-\frac{p+1}{2}}[x_{pp}-X_{21}X_{11}^{-1}X_{12}]^{\alpha-\frac{p+1}{2}}\eqno(1)
$$
$$|I-X|^{\beta-\frac{p+1}{2}}=|I-X_{11}|^{\beta-\frac{p+1}{2}}[(1-x_{pp})-X_{21}(I-X_{11})^{-1}X_{12}]^{\beta-\frac{p+1}{2}}
\eqno(2)
$$From (1) we have $x_{pp}>X_{21}X_{11}^{-1}X_{12}$ and from (2) we have $x_{pp}<1-X_{21}(I-X_{11})^{-1}X_{12}$. That is, $X_{21}X_{11}^{-1}X_{12}<x_{pp}<1-X_{21}(I-X_{11})^{-1}X_{12}$. Let $y=x_{pp}-X_{21}X_{11}^{-1}X_{12}\Rightarrow {\rm d}y={\rm d}x_{pp}$ for fixed $X_{21},X_{11}$. Also, $0<y<b$ where
\begin{align}
b&=1-X_{21}X_{11}^{-1}X_{12}-X_{21}(I-X_{11})^{-1}X_{12}\nonumber\\
&=1-X_{21}X_{11}^{-\frac{1}{2}}(I-X_{11})^{-\frac{1}{2}}(I-X_{11})^{-\frac{1}{2}}X_{11}^{-\frac{1}{2}}X_{12}\nonumber\\
&=1-WW',~W=X_{21}X_{11}^{-\frac{1}{2}}(I-X_{11})^{-\frac{1}{2}}.\nonumber
\end{align}The second factor on the right in (2) becomes
$$[b-y]^{\beta-\frac{p+1}{2}}=b^{\beta-\frac{p+1}{2}}[1-\frac{y}{b}]^{\beta-\frac{p+1}{2}}.
$$Put $u=\frac{y}{b}$ for fixed $b$. Then the factors containing $u$ and $b$ become $b^{\alpha+\beta-(p+1)+1}u^{\alpha-\frac{p+1}{2}}(1-u)^{\beta-\frac{p+1}{2}}$. Integral over $u$ gives
$$\int_0^1u^{\alpha-\frac{p+1}{2}}(1-u)^{\beta-\frac{p+1}{2}}{\rm d}u=\frac{\Gamma(\alpha-\frac{p-1}{2})\Gamma(\beta-\frac{p-1}{2})}{\Gamma(\alpha+\beta-(p-1))},
$$for $~\Re(\alpha)>\frac{p-1}{2},
~\Re(\beta)>\frac{p-1}{2}.$
 Let $W=X_{21}X_{11}^{-\frac{1}{2}}(I-X_{11})^{-\frac{1}{2}}$ for fixed $X_{11}$. Then ${\rm d}X_{21}=|X_{11}|^{\frac{1}{2}}|I-X_{11}|^{\frac{1}{2}}{\rm d}W$ from Theorem 1.18 of [4], where $X_{11}$ is $(p-1)\times (p-1)$. Put $v=WW'$ and integrate out over the Stiefel manifold by using Theorem 2.16 and Remark 2.13 of [4]. Then we have
$${\rm d}W=\frac{\pi^{\frac{p-1}{2}}}{\Gamma(\frac{p-1}{2})}v^{\frac{p-1}{2}-1}{\rm d}v.
$$Now the integral over $b$ becomes
\begin{align}
\int b^{\alpha+\beta-p}{\rm d}X_{21}&=\int_0^1v^{\frac{p-1}{2}-1}(1-v)^{\alpha+\beta-p}{\rm d}v\nonumber\\
&=\frac{\Gamma(\frac{p-1}{2})\Gamma(\alpha+\beta-(p-1))}{\Gamma(\alpha+\beta-\frac{p-1}{2})},~\Re(\alpha+\beta)>p-1.\nonumber
\end{align}Now, multiplying all the factors together we have
$$|X_{11}^{(1)}|^{\alpha+\frac{1}{2}-\frac{p+1}{2}}|I-X_{11}^{(1)}|^{\beta+\frac{1}{2}-\frac{p+1}{2}}\pi^{\frac{p-1}{2}}
\frac{\Gamma(\alpha-\frac{p-1}{2})\Gamma(\beta-\frac{p-1}{2})}{\Gamma(\alpha+\beta-\frac{p-1}{2})}
$$for $\Re(\alpha)>\frac{p-1}{2},~\Re(\beta)>\frac{p-1}{2}$. Here $X_{11}^{(1)}$ indicates the $(p-1)\times(p-1)$ leading submatrix at the end of the first set of operations. At the end of the second set of operations we will denote the $(p-2)\times (p-2)$ leading submatrix by $X_{11}^{(2)}$, and so on. The second step of operations starts by separating $x_{p-1,p-1}$ and writing

$$|X_{11}^{(1)}|=|X_{11}^{(2)}|[x_{p-1,p-1}-X_{21}^{(2)}[X_{11}^{(2)}]^{-1}X_{12}^{(2)}]
$$where $X_{21}^{(2)}$ is $1\times (p-2)$. Now, proceed as in the first sequence of steps to obtain the final factors of the following form:
$$|X_{11}^{(2)}|^{\alpha+1-\frac{p+1}{2}}|I-X_{11}^{(2)}|^{\beta+1-\frac{p+1}{2}}\pi^{\frac{p-2}{2}}
\frac{\Gamma(\alpha-\frac{p-2}{2})\Gamma(\beta-\frac{p-2}{2})}{\Gamma(\alpha+\beta-\frac{p-2}{2})}
$$for $\Re(\alpha)>\frac{p-2}{2},~\Re(\beta)>\frac{p-2}{2}$. Proceeding like this the exponent of $\pi$ at the end will be of the form
$$\frac{p-1}{2}+\frac{p-2}{2}+...+\frac{1}{2}=\frac{p(p-1)}{4}.
$$The gamma product will be of the form
$$\frac{\Gamma(\alpha-\frac{p-1}{2})\Gamma(\alpha-\frac{p-2}{2})...\Gamma(\alpha)\Gamma(\beta-\frac{p-1}{2})...\Gamma(\beta)}
{\Gamma(\alpha+\beta-\frac{p-1}{2})...\Gamma(\alpha+\beta)}.
$$These gamma products, together with $\pi^{\frac{p(p-1)}{4}}$ can be written as $\frac{\Gamma_p(\alpha)\Gamma_p(\beta)}{\Gamma_p(\alpha+\beta)}=B_p(\alpha,\beta)$ and hence the result. Thus, it is possible to evaluate the type-1 real matrix-variate beta integral directly to obtain the beta function in the real matrix variate case.
\vskip.2cm A similar procedure can yield the real matrix-variate beta function from the type-2 real matrix-variate beta integral of the form
$$\int_{X>O}|X|^{\alpha-\frac{p+1}{2}}|I+X|^{-(\alpha+\beta)}{\rm d}X
$$for $X=X'>O$ and $p\times p,~\Re(\alpha)>\frac{p-1}{2},~\Re(\beta)>\frac{p-1}{2}$. The procedure for the evaluation will be parallel.

\vskip.3cm\noindent{\bf 2.2.\hskip.3cm Evaluation of matrix-variate type-1 beta integral in the complex case}
\vskip.3cm The integral representation for $B_p(\alpha,\beta)$ in the complex case is the following:

$$\int_{O<\tilde{X}<I}|{\rm det}(\tilde{X})|^{\alpha-p}|{\rm det}(I-\tilde{X})|^{\beta-p}{\rm d}\tilde{X}=\tilde{B}_p(\alpha,\beta)
$$for $\Re(\alpha)>p-1,~\Re(\beta)>p-1$ where ${\rm det}(\cdot)$ denotes the determinant of $(\dot)$ and $|{\rm det}(\cdot)|$ denotes the absolute value of the determinant of $(\cdot)$. Here $\tilde{X}=(\tilde{x}_{ij})$ is a $p\times p$ hermitian positive definite matrix and hence all the diagonal elements are real and positive. As in the real case, let us separate $x_{pp}$ by partitioning:

$$\tilde{X}=\left[\begin{matrix}\tilde{X}_{11}&\tilde{X}_{12}\\
\tilde{X}_{21}&\tilde{X}_{22}\end{matrix}\right]\mbox{  as well as  }I-\tilde{X}=\left[\begin{matrix}I-\tilde{X}_{11}&-\tilde{X}_{12}\\
-\tilde{X}_{21}&I-\tilde{X}_{22}\end{matrix}\right].
$$Then the absolute value of the determinants are of the form:

$$|{\rm det}(\tilde{X})|^{\alpha-p}=|{\rm det}(\tilde{X}_{11})|^{\alpha-p}|x_{pp}-\tilde{X}_{21}\tilde{X}_{11}^{-1}\tilde{X}_{12}^{*}|^{\alpha-p}\eqno(a)
$$where * indicates conjugate transpose, and
$$|{\rm det}(I-\tilde{X})|^{\beta-p}=|{\rm det}(I-\tilde{X}_{11})|^{\beta-p}|(1-x_{pp})-\tilde{X}_{21}(I-\tilde{X}_{11})^{-1}\tilde{X}_{12}^{*}|^{\beta-p}.\eqno(b)
$$Note that when $\tilde{X}$ and $I-\tilde{X}$ are hermitian positive definite then $\tilde{X}_{11}^{-1}$ and $(I-\tilde{X}_{11})^{-1}$ are also hermitian positive definite. Further, the hermitian forms $\tilde{X}_{21}\tilde{X}_{11}^{-1}\tilde{X}_{12}^{*}$ and $\tilde{X}_{21}(I-\tilde{X}_{11})^{-1}\tilde{X}_{12}^{*}$ remain real and positive. From (a) and (b) it follows that

$$\tilde{X}_{21}\tilde{X}_{11}^{-1}\tilde{X}_{12}^{*}<x_{pp}<1-\tilde{X}_{21}(I-\tilde{X}_{11})^{-1}\tilde{X}_{12}^{*}.
$$Since hermitian forms are real, the lower and upper bounds of $x_{pp}$ are real. Let
$$\tilde{W}=\tilde{X}_{21}\tilde{X}_{11}^{-\frac{1}{2}}(I-\tilde{X}_{11})^{-\frac{1}{2}}
$$for fixed $\tilde{X}_{11}$. Then
$${\rm d}\tilde{X}_{21}=|{\rm det}(\tilde{X}_{11})|^{-1}|{\rm det}(I-\tilde{X}_{11})|^{-1}{\rm d}\tilde{W}
$$and $|{\rm det}(\tilde{X})|^{\alpha-p}, |{\rm det}(I-\tilde{X}_{11})|^{\beta-p}$ change to $|{\rm det}(\tilde{X}_{11})|^{\alpha+1-p}, |{\rm det}(I-\tilde{X}_{11})|^{\beta+1-p}$ respectively. Then we can write
\begin{align}
|(1-x_{pp})&-\tilde{X}_{21}\tilde{X}_{11}^{-1}\tilde{X}_{12}^{*}-\tilde{X}_{21}(I-\tilde{X}_{11})^{-1}\tilde{X}_{12}^{*}|^{\beta-p}\nonumber\\
&=(b-y)^{\beta-p}=b^{\beta-p}[1-\frac{y}{b}]^{\beta-p}.\nonumber
\end{align}Put $u=\frac{y}{b}$. Then the factors containing $u$ and $b$ will be of the form $u^{\alpha-p}(1-u)^{\beta-p}b^{\alpha+\beta-2p+1}$ and the integral over $u$ gives

$$\int_0^1u^{\alpha-p}(1-u)^{\beta-p}{\rm d}u=\frac{\Gamma(\alpha-(p-1))\Gamma(\beta-(p-1))}{\Gamma(\alpha+\beta-2(p-1))},
$$for $\Re(\alpha)>p-1,\Re(\beta)>p-1$.
 Let $v=\tilde{W}\tilde{W}^{*}$ and integrate out over the Stiefel manifold by using Corollaries 4.5.2 and 4.5.3 of [4]. Then

$${\rm d}\tilde{W}=\frac{\pi^{p-1}}{\Gamma(p-1)}v^{(p-1)-1}{\rm d}v.
$$The integral over $b$ gives
\begin{align}
\int b^{\alpha+\beta-2p+1}{\rm d}\tilde{X}_{21}&=\int_0^1v^{(p-1)-1}(1-v)^{\alpha+\beta-2p+1}{\rm d}v\nonumber\\
&=\frac{\Gamma(p-1)\Gamma(\alpha+\beta-2p+2)}{\Gamma(\alpha+\beta-p+1)},\nonumber
\end{align} for $\Re(\alpha)>p-1,\Re(\beta)>p-1$.
 Now, taking the product of all factors we have

$$|{\rm det}(\tilde{X}_{11})|^{\alpha+1-p}|{\rm det}(I-\tilde{X}_{11})|^{\beta+1-p}\pi^{p-1}\frac{\Gamma(\alpha-p+1)\Gamma(\beta-p+1)}{\Gamma(\alpha+\beta-p+1)}
$$for $\Re(\alpha)>p-1,\Re(\beta)>p-1$. Separate $x_{p-1,p-1}$ from $\tilde{X}_{11}$ and $I-\tilde{X}_{11}$ and continue the process. Then at the end, the exponent of $\pi$ will be $(p-1)+(p-2)+...+1=\frac{p(p-1)}{2}$ and the gamma product will be
$$\frac{\Gamma(\alpha-(p-1))\Gamma(\alpha-(p-2))...\Gamma(\alpha)\Gamma(\beta-(p-1))...\Gamma(\beta)}
{\Gamma(\alpha+\beta-(p-1))...\Gamma(\alpha+\beta)}.
$$These factors, together with $\pi^{\frac{p(p-1)}{2}}$ give
$$\frac{\tilde{\Gamma_p}(\alpha)\tilde{\Gamma_p}(\beta)}{\tilde{\Gamma_p}(\alpha+\beta)}=\tilde{B_p}(\alpha,\beta),
\Re(\alpha)>p-1,\Re(\beta)>p-1.
$$The procedure for evaluating a type-2 matrix-variate beta integral by the method of partitioning is parallel and hence it will not be detailed here.

\vskip.3cm\noindent{\bf 3.\hskip.3cm General Partitions}
\vskip.3cm In section 2 we have considered integrating one variable at a time by suitably partitioning the matrices. Is it possible to have a general partitioning and integrate a block of variables at a time, rather than integrating out individual variables? Let us consider the real matrix-variate gamma integral first. Let

$$X=\left[\begin{matrix}X_{11}&X_{12}\\
X_{21}&X_{22}\end{matrix}\right],~X_{11}\mbox{ is }p_1\times p_1\mbox{ and }X_{22}\mbox{ is }p_2\times p_2
$$so that $X_{12}$ is $p_1\times p_2$ and $X_{21}=X_{12}'$ and $p_1+p_2=p$. Without loss of generality, let us assume that $p_1\ge p_2$. Then the determinant can be partitioned as follows:
\begin{align}
|X|^{\alpha-\frac{p+1}{2}}&=|X_{11}|^{\alpha-\frac{p+1}{2}}|X_{22}-X_{21}X_{11}^{-1}X_{12}|^{\alpha-\frac{p+1}{2}}\nonumber\\
&=|X_{11}|^{\alpha-\frac{p+1}{2}}|X_{22}|^{\alpha-\frac{p+1}{2}}|I-X_{22}^{-\frac{1}{2}}X_{21}X_{11}^{-1}X_{12}
X_{22}^{-\frac{1}{2}}|^{\alpha-\frac{p+1}{2}}.\nonumber
\end{align}Put
$$Y=X_{22}^{-\frac{1}{2}}X_{21}X_{11}^{-\frac{1}{2}}\Rightarrow {\rm d}Y=|X_{22}|^{-\frac{p_1}{2}}|X_{11}|^{-\frac{p_2}{2}}{\rm d}X_{21}
$$for fixed $X_{11}$ and $X_{22}.$
$$|X|^{\alpha-\frac{p+1}{2}}=|X_{11}|^{\alpha+\frac{p_2}{2}-\frac{p+1}{2}}|X_{22}|^{\alpha+\frac{p_1}{2}-\frac{p+1}{2}}
|I-YY'|^{\alpha-\frac{p+1}{2}}.
$$The Jacobian above is available from Theorem 1.18 of [4]. Let $S=YY'$. Then integrating out over the Stiefel manifold we have

$${\rm d}Y=\frac{\pi^{\frac{p_1p_2}{2}}}{\Gamma_{p_2}(\frac{p_1}{2})}|S|^{\frac{p_1}{2}-\frac{p_2+1}{2}}{\rm d}S,
$$see Theorem 2.16 and Remark 2.13 of [4]. Now, integral over $S$ gives

$$\int_{O<S<I}|S|^{\frac{p_1}{2}-\frac{P_2+1}{2}}|I-S|^{\alpha-\frac{p_1}{2}-\frac{p_2+1}{2}}{\rm d}S=\frac{\Gamma_{p_2}(\frac{p_1}{2})\Gamma_{p_2}(\alpha-\frac{p_1}{2})}{\Gamma_{p_2}(\alpha)},
$$for $\Re(\alpha)>\frac{p_1-1}{2}$.
 Collecting all the factors, we have

$$|X_{11}|^{\alpha-\frac{p_1+1}{2}}|X_{22}|^{\alpha-\frac{p_2+1}{2}}\pi^{\frac{p_1p_2}{2}}\frac{\Gamma_{p_2}(\alpha
-\frac{p_1}{2})}{\Gamma_{p_2}(\alpha)}.
$$From here one can also observe that the original determinant splits into functions of $X_{11}$ and $X_{22}$. This also shows that if we are considering a real matrix-variate gamma density then the diagonal blocks $X_{11}$ and $X_{22}$ are statistically independently distributed, where $X_{11}$ will have a $p_1$-variate gamma distribution and $X_{22}$ has a $p_2$-variate gamma distribution. Observe that ${\rm tr}(X)={\rm tr}(X_{11})+{\rm tr}(X_{22})$ and hence the integral over $X_{22}$ gives $\Gamma_{p_2}(\alpha)$ and the integral over $X_{11}$ gives $\Gamma_{p_1}(\alpha)$. Hence the total integral is available as

$$\Gamma_{p_1}(\alpha)\Gamma_{p_2}(\alpha)\pi^{\frac{p_1p_2}{2}}\frac{\Gamma_{p_2}(\alpha-\frac{p_1}{2})}
{\Gamma_{p_2}(\alpha)}
=\Gamma_p(\alpha)$$since
 $\pi^{\frac{p_1p_2}{2}}\Gamma_{p_1}(\alpha)\Gamma_{p_2}(\alpha-\frac{p_1}{2})=\Gamma_p(\alpha)$. \vskip.2cm Hence it is seen that instead of integrating out variables one at a time we could have also integrated out blocks of variables at a time and could have verified the result. Similar procedure works for real matrix-variate type-1 and type-2 beta, and matrix-variate gamma, type-1 and type-2 beta in the complex domain also.

 \vskip.3cm\noindent{\bf 3.1.\hskip.3cm Methods avoiding integration over the Stiefel manifold}

 \vskip.3cm The general method of partitioning described above involves the integration over the Stiefel manifold as an intermediate step. We will consider another procedure which will avoid integration over Stiefel manifold. Let us consider the real gamma case first. Again, we start with the decomposition
 $$|X|^{\alpha-\frac{p+1}{2}}=|X_{11}|^{\alpha-\frac{p+1}{2}}|X_{22}-X_{21}X_{11}^{-1}X_{12}|^{\alpha-\frac{p+1}{2}}.\eqno(3.1)
 $$Instead of integrating out $X_{21}$ or $X_{12}$ let us integrate out $X_{22}$. Let $X_{11}$ be $p_1\times p_1$ and $X_{22}$ be $p_2\times p_2$ with $p_1+p_2=p$. In the above partitioning we require that $X_{11}$ be nonsingular. But when $X$ is positive definite, both $X_{11}$ and $X_{22}$ will be positive definite, thereby nonsingular also. From the second factor in (3.1), $X_{22}>X_{21}X_{11}^{-1}X_{12}$ from $X_{22}-X_{21}X_{11}^{-1}X_{12}$ being positive definite. We will try to integrate out $X_{22}$ first. Let $U=X_{22}-X_{21}X_{11}^{-1}X_{12}$ so that ${\rm d}U={\rm d}X_{22}$ for fixed $X_{11}$ and $X_{12}$. Since ${\rm tr}(X)={\rm tr}(X_{11})+{\rm tr}(X_{22})$ we have
 $${\rm e}^{-{\rm tr}(X_{22})}={\rm e}^{-{\rm tr}(U)-{\rm tr}(X_{21}X_{11}^{-1}X_{12})}.
 $$Integrating out $U$ we have
 $$\int_{U>O}|U|^{\alpha-\frac{p+1}{2}}{\rm e}^{-{\rm tr}(U)}{\rm d}U=\Gamma_{p_2}(\alpha-\frac{p_1}{2}),~\Re(\alpha)>\frac{p-1}{2}
 $$since $\alpha-\frac{p+1}{2}=\alpha-\frac{p_1}{2}-\frac{p_2+1}{2}$. Let
 $$Y=X_{21}X_{11}^{-\frac{1}{2}}\Rightarrow {\rm d}Y=|X_{11}|^{-\frac{p_2}{2}}{\rm d}X_{21}
 $$for fixed $X_{11}$. Then
 $$\int_{X_{21}}{\rm e}^{-{\rm tr}(X_{21}X_{11}^{-1}X_{12})}{\rm d}X_{21}=|X_{11}|^{\frac{p_2}{2}}\int_{Y}{\rm e}^{-{\rm tr}(YY')}{\rm d}Y.
 $$But ${\rm tr}(YY')$ is the sum of squares of the $p_1p_2$ elements in $Y$ and each integral is of the form $\int_{-\infty}^{\infty}{\rm e}^{-z^2}{\rm d}z=\sqrt{\pi}$. Hence
 $$\int_{Y}{\rm e}^{-{\rm tr}(YY')}{\rm d}Y=\pi^{\frac{p_1p_2}{2}}.
 $$Now we can integrate out $X_{11}$.
 \begin{align}
 \int_{X_{11}>O}|X_{11}|^{\alpha+\frac{p_2}{2}-\frac{p+1}{2}}&{\rm e}^{-{\rm tr}(X_{11})}{\rm d}X_{11}\nonumber\\
 &=\int_{X_{11}>O}|X_{11}|^{\alpha-\frac{p_1+1}{2}}{\rm e}^{-{\rm tr}(X_{11})}{\rm d}X_{11}\nonumber\\
 &=\Gamma_{p_1}(\alpha).\nonumber
 \end{align}Hence we have the following factors:
 $$\pi^{\frac{p_1p_2}{2}}\Gamma_{p_2}(\alpha-\frac{p_1}{2})\Gamma_{p_1}(\alpha)=\Gamma_p(\alpha)
 $$since
 $$\frac{p_1(p_1-1)}{4}+\frac{p_2(p_2-1)}{4}+\frac{p_1p_2}{2}=\frac{p(p-1)}{4},~p=p_1+p_2
 $$and
 \begin{align}
 \Gamma_{p_1}(\alpha)\Gamma_{p_2}(\alpha-\frac{p_1}{2})&=\Gamma(\alpha)\Gamma(\alpha-\frac{1}{2})
 ...\Gamma(\alpha-\frac{p_1-1}{2})\Gamma_{p_2}(\alpha-\frac{p_1}{2})\nonumber\\
 &=\Gamma(\alpha)...\Gamma(\alpha-\frac{p_1+p_2-1}{2}).\nonumber
 \end{align}Hence the result. In this procedure we did not have to go through integration over the Stiefel manifold and we did not have to assume that $p_1\ge p_2$. We could have integrated out $X_{11}$ first if needed. In this case, expand
 $$|X|^{\alpha-\frac{p+1}{2}}=|X_{22}|^{\alpha-\frac{p+1}{2}}|X_{11}-X_{12}X_{22}^{-1}X_{21}|^{\alpha-\frac{p+1}{2}}.
 $$Then proceed as before by integrating out $X_{11}$ first. Then we end up with
 $$\pi^{\frac{p_1p_2}{2}}\Gamma_{p_1}(\alpha-\frac{p_2}{2})\Gamma_{p_2}(\alpha)=\Gamma_p(\alpha), p=p_1+p_2.
 $$
 \vskip.3cm\noindent{\bf Note:}\hskip.3cm If we are considering a real matrix-variate gamma density, such as the Wishart density, then from the above procedure observe that after integrating out $X_{22}$ the only factor containing $X_{21}$ is the exponential function, which has the structure of a matrix-variate Gaussian density. Hence for given $X_{11}$, $X_{21}$ is matrix-variate Gaussian distributed. Similarly, for given $X_{22}$, $X_{12}$ is matrix-variate Gaussian distributed. Further, the diagonal blocks $X_{11}$ and $X_{22}$ are independently distributed.
 \vskip.3cm The same procedure as above goes through for the evaluation of gamma integrals in the complex domain also. Since the steps are parallel they will not be detailed here.

\vskip.3cm\noindent{\bf Acknowledgement}
\vskip.3cm The author would like to thank the Department of Science and Technology, Government of India, for the financial assistance for this work under project number SR/S4/MS:287/05 and the Centre for Mathematical Sciences for the facilities.

\vskip.3cm\noindent\begin{center}
References
\end{center}
\vskip.3cm\noindent [1]~~T.W. Anderson,  An Introduction to Multivariate Statistical Analysis, Wiley, New York, 1971.
\vskip.2cm\noindent [2]~~A.M. Kshirsagar,  Multivariate Analysis, Marcel Dekker, New York, 1972.
\vskip.2cm\noindent [3]~~M.S. Srivastava and C.G., An Introduction to Multivariate Statistics, North Holland, New York, 1979.
\vskip.2cm\noindent [4]~~A.M. Mathai, Jacobians of Matrix Transformations and Functions of Matrix Argument, World Scientific Publishing, New York, 1997.

\end{document}